\documentclass[12pt]{amsart}
\usepackage{epsfig,amsfonts,amssymb,amsthm}
\usepackage{amsmath}
\usepackage{fullpage}

\usepackage{amssymb}
\usepackage{epsfig,amsfonts,amssymb,amsthm}

\begin{document}




\title{Bandgaps in two-dimensional high-contrast periodic elastic beam lattice materials}


\author{I.V.
Kamotski}
\address{Department of Mathematics, University College London, Gower Street, London WC1E 6BT, UK}

\author{V. P. Smyshlyaev}
\address{Department of Mathematics, University College London, Gower Street, London WC1E 6BT, UK}
  \email{v.smyshlyaev@ucl.ac.uk}

\dedicatory{Dedicated to Professor Norman A. Fleck on occasion of his 60-th birthday}

\begin{abstract}
We consider elastic waves in a two-dimensional periodic lattice network of Timoshenko-type beams. 
We show that for general configurations involving certain highly-contrasting components a high-contrast  modification of the homogenization theory is capable of accounting for bandgaps, explicitly relating 
those to low resonant frequencies of the ``soft'' components. 
An explicit example of a square-periodic network of beams with a single isolated resonant beam within a periodicity cell is considered in detail. 

\end{abstract}






\maketitle

\section{Introduction}
\label{}

Macroscopic properties of lattice materials, and their static and dynamic responses have been of a
considerable recent interest. In particular, Phani, Woodhouse and Fleck \cite{PhFl06} studied wave propagation for certain
lattice topologies and associated phenomena of frequency bandgaps and spatial filtering (wave
directionality) by adopting Floquet-Bloch's principles for lattices modeled as a network of 
Timoshenko beams
using the finite element method. Long-wavelength asymptotes to the resulting dispersion curves
were found to be in good agreement with those based on homogenization (effective medium) theories,
however the range of validity of the latter was found to be restricted to low frequencies i.e.
below the frequencies where the bandgaps may be observed.

In the present work we argue however that for periodic elastic beam configurations involving certain 
{\it highly-contrasting} components,  
a high-contrast modification of the homogenization theory is capable of accounting both for
bandgaps and for certain forms of polarization filtering, with the former explicitly related to (low) resonance
frequencies of the ``soft'' components. As a result of the high contrast, and as opposed to the 
classical homogenization, under a naturally chosen micro-resonant scaling the asymptotic description 
of the wave processes remains intrinsically two-scale, which in fact gives rise to the above effects. 
We emphasize that the mathematical approach adopted here is rigorous in the sense that 
it does not rely on any further assumptions, and provides approximations for entities of interest 
for the underlying Floquet-Bloch waves and in particular for the bandgaps with 
a controllably small error for a sufficiently high value of the parameter of the contrast. 

Without attempting here a comprehensive review of the relevant literature,
some related ideas for two-phase high-contrast and generally highly-anisotropic periodic elastic composites were discussed in \cite{VPS09} with some examples 
of both frequency gaps and of long-wavelength ``directional localization'', following mathematical
ideas of \cite{Zhikov00,Zhikov05}. Related developments specifically for high-contrast linear elasticity 
include those by \cite{bellieud}, \cite{ZhP13}, \cite{Cooper14}, and for high-contrast graphs but outside elasticity 
by \cite{Cheredn18}. 
The derived macroscopic equations in the bandgap regime allow interpretation in terms of negative 
effective density (or sign-indefinite anisotropic density in case of ``weak'' gaps) for some frequency ranges. 
Some related general ideas are found in \cite{Willis85}, 
and in a specific context of high-frequency periodic homogenization 
 this was probably 
first observed by \cite{AurBonn85,Aur94} and developed further by \cite{BF04} and \cite{Avila}, among others. 
Some related general mathematical issues for two-scale homogenization
of systems of partial differential equations with partially degenerating periodic coefficients were
recently studied by us in \cite{KSAA18}. 

From the mathematical perspective, the case of high-contrast 
 Timoshenko-type elastic beam networks appears of additional interest (and indeed of a non-trivial additional challenge to us)  for the non-trivial effect of 
microscopic rotations on the two-scale limit asymptotic behavior. 
In particular, the three microscopic degrees of freedom (two displacements and a rotation) generate, 
for a given frequency, {\it up to three} propagating Bloch modes. However, in the chosen (two-scale) high-contrast 
asymptotic regime, only {\it up to two} modes can propagate, as in the conventional two-dimensional (2-D)
elasticity: the rotational degree of freedom remains purely microscopic, but nevertheless still affects 
the macroscopic part via the two-scale coupling. 
The latter is due to the fact that, for a periodically connected stiff component, the homogenized 
tensor turns out to be conventional 2-D elastic (i.e. with no macroscopic rotational degree of freedom): this is not obvious a priori, and we establish this as a by-product of our 
analysis. 
A further related feature of the limit two-scale model is that, due to the coupling between the two above 
macroscopic modes, in principle 
the microscopic resonances may not necessarily lead to bandgaps (in contrast to the prototype scalar 
models of \cite{Zhikov05}). We show that the gaps nevertheless do appear, at least for 
particular configurations.

The issues of filtering properties  
for various models of elastic lattice structures including band gaps have also been intensively studied before, 
see for example \cite{Mead96} and further references therein, and 
\cite{MM03}, although to the best of our knowledge not specifically for the 
high-contrast ``micro-resonant'' scaling for which, as we argue in this work, the lattice structures with highly-contrasting components 
provide distinct scenarios for such effects and for their understanding in terms of the underlying
lattice resonances. On the other hand, effect of the resonances on the bandgaps has been studied e.g. by  
\cite{Ph13} for a simple model of a single Timoshenko beam with attached periodic masses and
resonators. 

The structure of the paper is as follows. In Section 2 we describe a model for wave propagation in a 
two-dimensional periodic elastic Timoshenko beam network, which is essentially the same as the 
underlying initial model of rigid-jointed network of beams of \cite{PhFl06}, but without any further 
approximations for the fields on the beams and 
 with some high-contrast elements. We then explain how 
the resonance effect for the soft part 
leads to a natural scaling for high contrast vs small periodicity. Section 3 
executes two-scale asymptotic analysis of the emerging problem, and derives a two-scale limit problem. 
The limit problem is analyzed in Section 4, which discusses how it reveals the effect of full and partial (weak) band gaps and their relation to the underlying resonances. Section 5 considers an explicit example of a 
square beam network with isolated resonant beams, where the entities of interest are  
evaluated analytically, and as a result the existence of both full and ``partial'' bandgaps is established 
via a further qualitative argument. The Appendix explicitly calculates the homogenized elasticity tensor for the 
stiff component in the above example. 

\section{Periodic beam networks with a high contrast} 

Similarly to \cite{PhFl06}, 
we consider a two-dimensional periodic rigid-jointed network of beams with no pre-stress akin to the one displayed on Fig. 1. 
\begin{figure}
	\centering
		\includegraphics[scale=0.5]{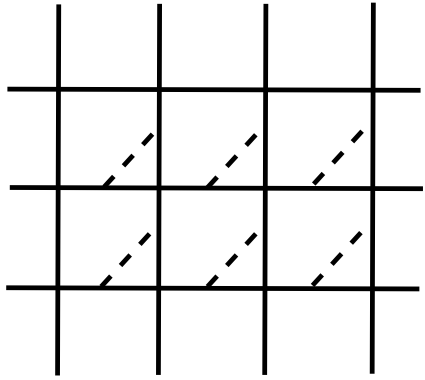}
		\caption{Geometric configuration: a periodic network of beams} 
		\label{fig1} 
\end{figure} 
Each beam of a length $L$ is modeled as a Timoshenko beam, see e.g. \cite{WJ87}, with any of its center-line  material 
points $x$, $-L/2\leq x\leq L/2$, having three degrees of freedom (Fig. \ref{fig2}): longitudinal displacement $u(x)$ along the 
beam, displacement $v(x)$ in the transverse direction $y$, and the total rotation $\theta(x)$ 
of the material normal to the undeformed beam 
about the $z$-axis. 
(The latter is a combined effect of the rotation $dv/dx$ of the beam and of an additional rotation due to the 
shear within the beam.)

\begin{figure}
	\centering
		\includegraphics[scale=0.27]{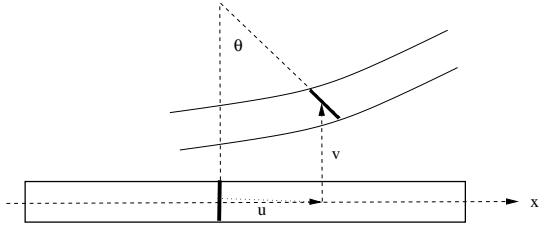}
		\caption{Timoshenko beam} 
		\label{fig2} 
\end{figure} 

The resulting kinetic and potential energies per unit thickness of the beam in the $z$-direction, are: 
\begin{equation}
T_b\,=\,\frac{1}{2}\int_{-L/2}^{L/2}\left[\,\rho\dot{u}^2\,+\rho\dot{v}^2\,+\,\rho I\dot{\theta}^2\right]\,dx,  
\label{Tkin}
\end{equation} 
\begin{equation}
U_b\,=\,\frac{1}{2}\int_{-L/2}^{L/2}
\left[\,E\left(\frac{du}{dx}\right)^2\,+\,E I\left(\frac{d\theta}{dx}\right)^2
\,+\,G  \left(\frac{d v}{dx}-\theta\right)^2
\right]\,dx.  
\label{Upot}
\end{equation} 
Here the positive parameters describe physical and geometrical properties of the beam: 
$\rho$ is the density of the beam per unit length, $I$ is the second moment of area of the beam divided by its 
thickness, 
and $E$ and $G$ characterize mechanical stiffness of the beam in extension and shear, see e.g. 
\cite{WJ87} and \cite{PhFl06} for more details. These parameters may vary from one beam to another within a periodicity 
cell, however they are assumed to be periodically replicated from one cell to another. 
Notice that $U_b=0$ if and only if $\theta=c_1$, $u=c_2$ and $v=c_3+c_1y$ with constant 
$c_1$, $c_2$ and $c_3$, which corresponds to rigid body motions. 
The imposed below 
following \cite{PhFl06} 
condition of rigid joints implies the continuity of the relevant components of displacement as 
well as the continuity of the rotation $\theta$. 

The total kinetic and potential energies within a finite volume $V$ of the network are the sums over 
all beams within $V$, 
\[
T\,=\,\sum_j T_{b_j}, \ \ \ \ U\,=\,\sum_j U_{b_j}, 
\]
and the equations of motion, in the absence of external forces, are derived by applying Hamilton's variational principle, 
\begin{equation}
\delta \int (T\,-\,U)\,dt\,=\,0. 
\label{hamilton}
\end{equation} 
The latter, via (\ref{Tkin}) and (\ref{Upot}), can be conveniently written in the following weak form: 
\[
\sum_j\int_{-L_j/2}^{L_j/2}
\left[
\,\rho_j\ddot{u}_j\tilde u_j\,+\rho_j\ddot{v}_j\tilde v_j\,+\,\rho_j I_j\ddot{\theta}_j\tilde\theta_j\,\,\,\,+ 
\ \ \ \ \ \ \ \ \ \ \ \ \ \  
\right. 
\]
\begin{equation}
\left. 
\,E_j \frac{du_j}{dx_j}\frac{d\tilde u_j}{dx_j}\,+\,E_j I_j\frac{d\theta_j}{dx_j}\frac{d\tilde\theta_j}{dx_j}
\,+\,G_j  \left(\frac{d v_j}{dx_j}-\theta_j\right)\left(\frac{d \tilde v_j}{dx_j}-\tilde\theta_j\right)
\right]\,dx_j\,\,=\,\,0.    
\label{eqmovar}
\end{equation} 
Here the integral identity is required to hold, for any time $t$, 
for all smooth test functions 
$(\tilde u_j,\tilde v_j,\tilde \theta_j)$  supported on the graph within $V$, 
which together with the sought solution $(u_j,v_j,\theta_j)$
satisfy the kinematic continuity 
conditions at the joints. Namely, at a joint $J$, for all the beams $b_j$ connecting to $J$ 
(denoted $j\in J$), $\theta_j$ are continuous and the displacements associated with $u_j$ and $v_j$ 
are continuous as well, i.e. 
\begin{equation}
\theta_j\,=\,\theta_l\,\,\mbox{ and }\,\, 
u_j\mbox{\boldmath{$\tau$}}^{(j)}\,+\,v_j\mathbf{n}^{(j)}\,=\,
u_l\mbox{\boldmath{$\tau$}}^{(l)}\,+\,v_l\mathbf{n}^{(l)},  \,\,\,
\mbox{ for all } j,l\in J, 
\label{contjunc}
\end{equation} 
where $\mbox{\boldmath{$\tau$}}^{(j)}$ is the unit vector along the beam $b_j$ in the direction of the increase of $u_j$ and $\mathbf{n}^{(j)}$
 is the unit vector in the normal direction of the positive $v_j$. 

We emphasize here that the chosen microscopic model is essentially identical to the underlying 
initial model of \cite{PhFl06} of rigid-jointed network of beams. In (\ref{contjunc}), together with the  
condition of continuity of the displacements (which is natural assuming the joint region is small enough), both require the continuity of 
the rotations. Notice that the latter continuity is implicitly enforced by \cite{PhFl06} 
via prescribing $\theta$ (as well as the total displacements) at the nodes only and then 
continuously 
\begin{figure}
	\centering
		\includegraphics[scale=0.27]{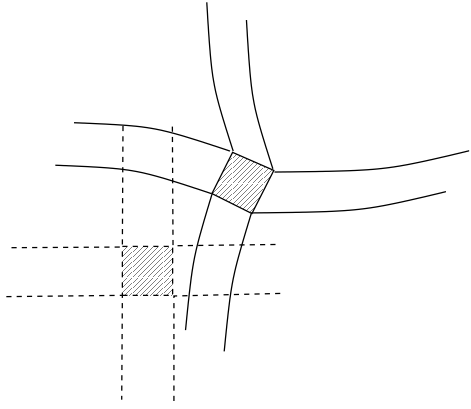}
		\caption{Rigid joint regions, shown as shaded squares for the undeformed (dashed lines) and the deformed 
		(solid lines) junction of four Timoshenko beams} 
		\label{fig2extra} 
\end{figure} 
approximating them between the nodes (i.e. on the beams) by particular ``shape functions'' (Appendix 
of \cite{PhFl06}), within the ``exact'' Hamilton's variational framework (\ref{hamilton}). Our approach remains 
exact in the sense that it does not involve any such further approximations. 
The continuity of the rotations appears natural, as it assumes that the rigid joint region is of a small size 
commensurable with the thickness $d$ of the beams, which is in turn assumed much smaller than the length 
of the lattice $L$. 
All the beams, 
see an illustration on Fig. \ref{fig2extra}, are assumed bonded to the joint along 
their undeformed normals and so, 
as only rigid motion (i.e. translations and rotations) of the joint is allowed, the normals of all the beams in the 
immediate vicinity of the joint have the same amount of rotation $\theta$. 
A systematic derivation of such a condition would require asymptotic analysis involving a small parameter of 
the size of the joint, cf e.g. \cite{MKM,Panas05} and/or variational approximations to Hamilton's principle 
(\ref{hamilton}) cf. e.g. \cite{Willis81,SCh00}, which is beyond the scope of the present work. 
Such a model of rigid joints is often postulated in the literature, e.g. in \cite{Lagnese93} p. 352, 
which leads to a mathematically well-posed problem on the resulting network. So we adopt this 
particular model as a 
part of our microscopic problem, although our approach can be applied to other models of beam 
networks as well.

For time-harmonic waves of an angular frequency $\omega$ which may propagate through such an 
infinite-periodic lattice, in particular for Floquet-Bloch waves, the equation of motion is 
obtained by formally adopting the time dependence throughout in (\ref{eqmovar}) through the factor 
$e^{i\omega t}$. Allowing the same notation for the time-harmonic part, i.e. assuming in 
(\ref{eqmovar}) $u_j(x,t)=e^{i\omega t} u_j(x)$, etc, transforms (\ref{eqmovar}) into a spectral problem 
\[
\sum_j\int_{-L_j/2}^{L_j/2}
\left[
\,E_j \frac{du_j}{dx_j}\frac{d\tilde u_j}{dx_j}\,+\,E_j I_j\frac{d\theta_j}{dx_j}\frac{d\tilde\theta_j}{dx_j}
\,+\,G_j  \left(\frac{d v_j}{dx_j}-\theta_j\right)\left(\frac{d \tilde v_j}{dx_j}-\tilde\theta_j\right)
\right]\,dx_j\,\,= 
\]
\begin{equation} 
\omega^2\,\sum_j\int_{-L_j/2}^{L_j/2}
\left[
\,\rho_j{u}_j\tilde u_j\,+\rho_j{v_j}\tilde v_j\,+\,\rho_j I_j {\theta_j}\tilde\theta_j\,
\right]\,dx_j. 
\label{eqspectralvar}
\end{equation} 
Formal integration by parts in (\ref{eqspectralvar}) gives equations of motion on each beam, as well 
as conditions on equilibrium of forces and torques at the joints. For example, for the 
longitudinal displacements $u_j$ the standard Helmholtz equation holds 
\[
\frac{d}{dx_j}\left(\,E_j\frac{du_j}{dx_j}\right)\,+\,\rho_j\omega^2u_j\,=\,0. 
\]
It is more convenient for us however to operate, up to certain point,  directly with the weak 
form (\ref{eqspectralvar}). 

Quasi-periodic solutions to (\ref{eqspectralvar}) are the Bloch waves. For example, for square lattices 
with period $l$, the quasi-periodicity condition with ``quasi-momentum'' $\mathbf{k}=(k_1,k_2)$ for $u$ reads 
$u(\mathbf{x}+\mathbf{m}l)=e^{i\mathbf{k}\cdot\mathbf{m}l}u(\mathbf{x})$ for any $\mathbf{m}=(m_1,m_2)$ 
with integer $m_1$ and $m_2$, with similar conditions for $v$ and $\theta$. Each $\mathbf{k}$ has 
associated discrete frequencies $\omega_j(\mathbf{k})$, $j=1,2,...$. Such $\omega$ which 
have no associated $\mathbf{k}$ (if any) are the band gap frequencies, i.e. the forbidden frequencies 
which cannot propagate through such a structure.

We are specifically interested in certain {\it high-contrast} lattice networks, namely such that a 
connected ``stiff'' periodic component e.g. a cubic lattice displayed by solid lines on Fig. \ref{fig1}, 
contains additional possibly disconnected periodic ``soft'' beam elements like those inclined beams 
displayed on Fig. 1 by dashed lines. Introducing a small parameter of the contrast $\delta$, $0<\delta \ll 1$, 
we assume that all (or, in some further generalizations, possibly some, cf \cite{VPS09}) 
the pre-factors on the left hand side of (\ref{eqspectralvar}) for the soft phase are order-$\delta$ smaller 
than those for the stiff phase, i.e. $E_{{{\rm soft}}}/E_{{\rm stiff}}\sim \delta$, 
$E_{{\rm soft}}I_{{\rm soft}}/(E_{{\rm stiff}}I_{{\rm stiff}})\sim \delta$, 
$G_{{\rm soft}}/G_{{\rm stiff}}\sim \delta$. 

For a high contrast, i.e. for a small positive $\delta$, we are interested in time harmonic waves 
which may propagate through the lattice such that their frequency $\omega$ is commensurable with the
main eigen-frequency of the soft components with ``clamped'' end points (i.e. with $u=v=\theta=0$ for 
$x=\pm L_{\rm soft}/2$). Because of the high contrast, such ``resonant'' frequencies with respect to the soft phase 
will be perceived as low frequencies in the surrounding stiff matrix: while the corresponding wavelength 
in the soft components will be comparable to the length $L_{\rm soft}$ of the soft beam, it is much 
larger than the periodicity size in the stiff part. Had there been no ``resonant'' soft parts, such a low 
frequency regime would allow employing appropriate effective medium theories, in particular homogenization 
theory, approximating the beam network by an equivalent continuum medium at a macroscale of the order 
of the wavelength in the stiff phase. The latter continuum however displays no band gap or other effects 
of our particular interest, and a key point for us here is that introducing the soft ``resonators'' does 
on one hand allow to keep relevance of the (appropriately modified) homogenization theory and on the other hand 
does allow in particular to observe the band gap effects. 

Keeping, for simplicity,  the high contrast only in the stiffness but not in 
density\footnote{Note that, for analogous continuous problem, some models which 
involved high contrast not only in stiffness but also in density were considered by \cite{BKS08}}, 
i.e. in the pre-factors of (\ref{Upot}) but not (\ref{Tkin}), 
observe 
that for the wavelength $\lambda$, $\lambda\sim\omega^{-1}\left(E/\rho\right)^{1/2}$. So for the wavelengths 
$\lambda_{\rm stiff}$ and $\lambda_{\rm soft}$ in the soft and in the stiff phases respectively, 
$\lambda_{\rm soft}/\lambda_{\rm stiff}\sim \left(E_{\rm soft}/E_{\rm stiff}\right)^{1/2}\sim\delta^{1/2}$. 
Therefore, setting $\varepsilon\sim\delta^{1/2}$, and choosing for the ``order one'' macroscale 
comparable to 
$\lambda_{\rm stiff}$, we have 
\[
O(1)\,\sim\,\lambda_{\rm stiff}\,\sim\,\varepsilon^{-1}\lambda_{\rm soft}\,\sim\,\varepsilon^{-1}L, 
\]
i.e. for the beams lengths $L$ and therefore for the periodicity size $l$, 
\[
l\,\sim\,L\,\sim\,\varepsilon, \ \  \ \mbox{with the contrast} \ \ \delta\,\sim\, \varepsilon^2. 
\]

The above dimensional analysis suggests that at the microscale all the coefficients in (\ref{eqspectralvar}) 
for the soft phase have to be chosen of order one, while for the stiff phase the coefficients on the 
right hand side of (\ref{eqspectralvar}) are still of order one but on the left hand side $E_j$, $E_jI_j$ and $G_j$ are of order 
$\varepsilon^{-2}$. Rescaling to the macroscale, i.e. changing $x$ to $x/\varepsilon$ (and hence 
$d/dx_j$ to $\varepsilon d/dx_j$), and assuming for 
notational simplicity all the remaining order-one pre-factors to be either constants or unities and 
regarding $\lambda:=\omega^2$ as 
a spectral parameter results in 
\[
\int_{\Gamma_{\rm stiff}^\varepsilon}\varepsilon 
\left[
\gamma\frac{du}{dx}\frac{d\tilde u}{dx}\,+\,\eta\frac{d\theta}{dx}\frac{d\tilde\theta}{dx}
\,+\, \kappa\left(\frac{d v}{dx}-\varepsilon^{-1}\theta\right)
\left(\frac{d \tilde v}{dx}-\varepsilon^{-1}\tilde\theta\right)
\right]\,dx\,\,+\,  \ \ \ \ \ \ \ \ \ \ \ \ \ \ \ \ \ \ \ 
\]
\[
\int_{\Gamma_{\rm soft}^\varepsilon}
\varepsilon^3
\left[
\,\frac{du}{dx}\frac{d\tilde u}{dx}\,+\,\frac{d\theta}{dx}\frac{d\tilde\theta}{dx}
\,+\, \left(\frac{d v}{dx}-\varepsilon^{-1}\theta\right)
\left(\frac{d \tilde v}{dx}-\varepsilon^{-1}\tilde\theta\right)
\right]\,dx\,\,\,\,\,\,=  \ \ \ \ \ \ \ \ \
\]
\begin{equation}  
\lambda\,\left(\int_{\Gamma_{\rm stiff}^\varepsilon}+\int_{\Gamma_{\rm soft}^\varepsilon}\right) 
\varepsilon 
\left[
\,{u}\tilde u\,+\,{v}\tilde v\,+\, {\theta}\tilde\theta\,
\right]\,dx. 
\label{eqspvarmacro}
\end{equation} 
In (\ref{eqspvarmacro}), 
$\gamma$, $\eta$ and $\kappa$ are positive constants; 
$\Gamma_{\rm stiff}^\varepsilon$ and $\Gamma_{\rm soft}^\varepsilon$ denote  
infinite $\varepsilon$-periodic stiff and soft graphs respectively, i.e. 
$\Gamma_{\rm stiff}^\varepsilon:=\varepsilon \Gamma_{\rm stiff}$ and 
$\Gamma_{\rm soft}^\varepsilon:=\varepsilon \Gamma_{\rm soft}$ where $\Gamma_{\rm stiff}$ and 
$\Gamma_{\rm soft}$ are related reference graphs of (double-)period one. 
For example, for the square lattice-type beam network on Fig. 1,  
$\Gamma_{\rm stiff}$ and 
$\Gamma_{\rm soft}$ are periodically replicated stiff (solid) and soft (dashed) 
parts of the unit cell graph $\Gamma$ displayed on Fig. 3, and consisting of a stiff component 
denoted $\Gamma_1$ (the solid cross) and of the soft component $\Gamma_0$ (the dashed segment). 
\begin{figure}
	\centering
		\includegraphics[scale=0.32]{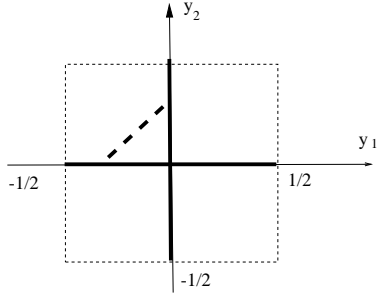}
		\caption{Unit cell graph $\Gamma$} 
		\label{fig3} 
\end{figure} 
Therefore $\Gamma_{\rm stiff}^\varepsilon$ is an $\varepsilon$-sized grid, consisting of all the point 
on the $(x,y)$-plane  
such that either $x=\varepsilon n$ with an integer $n$ and arbitrary $y$ or of $y=\varepsilon m$ with 
an integer $m$ and arbitrary $x$, and with appropriate $\Gamma_{\rm soft}^\varepsilon$. 
Further, $dx$ in (\ref{eqspvarmacro}) are line elements on $\Gamma_{\rm stiff}^\varepsilon$ and 
$\Gamma_{\rm soft}^\varepsilon$, $u$ are longitudinal displacements along $dx$, and $v$ are transverse 
displacements in the direction of $dx$ rotated 90$^{\circ}$ anti-clockwise; $\theta$ is rotation 
in the anti-clockwise direction, cf Fig. \ref{fig2}. 
Notice that in (\ref{eqspvarmacro}), whose both sides can be simultaneously pre-multiplied by any power of $\varepsilon$, 
we chose pre-factors $\varepsilon$ for the stiff parts (and so $\varepsilon^3$ for the soft parts): 
this is for ensuring that for order-one entities in the integrands' square brackets the integrals' values 
are also order-one. (As, within an order-one volume $V$, one has order $\varepsilon^{-2}$ many of 
$\varepsilon$-sized periodic cells with associated order-$\varepsilon$ line integrals.) 

\section{Two-scale homogenization for problem (\ref{eqspvarmacro})} 

We describe in this section a formal asymptotic procedure for solving the spectral problem 
(\ref{eqspvarmacro}) as 
$\varepsilon\to 0$. 
It follows general recipes of two-scale asymptotic expansions as extended to (high-contrast) 
periodic homogenization, cf e.g. \cite{BLP,BaPan,VPS09,KSJMS18}, as adapted in a natural way 
to periodic graph structures, so the 
derivation which we give below is relatively sketchy. 

One seeks a formal solution to the problem (\ref{eqspvarmacro}) in the form of a standard 
two-scale asymptotic ansatz as adapted to the $\varepsilon$-periodic network $\Gamma^\varepsilon$, 
the union of $\Gamma_{\rm stiff}^\varepsilon$ and $\Gamma_{\rm soft}^\varepsilon$. 
Namely, on $\Gamma^\varepsilon$, the displacement 
$\mathbf{u}^\varepsilon(\mathbf{x})=\left(u_1^\varepsilon(\mathbf{x}), u_2^\varepsilon(\mathbf{x})\right)$, 
$\mathbf{x}=(x_1,x_2)$, 
as well as the rotation $\theta^\varepsilon(\mathbf{x})$ are sought in the form of a two-scale 
asymptotic expansion 
\begin{eqnarray}
\mathbf{u}^\varepsilon(\mathbf{x})&\sim&\mathbf{u}^{(0)}(\mathbf{x},{x}/\varepsilon)+
\varepsilon\mathbf{u}^{(1)}(\mathbf{x},x/\varepsilon)+ 
\varepsilon^2\mathbf{u}^{(2)}(\mathbf{x},{x}/\varepsilon)+..., 
\label{u2sc}
\\ 
\theta^\varepsilon(\mathbf{x})&\sim&\theta^{(0)}(\mathbf{x},x/\varepsilon)+
\varepsilon\theta^{(1)}(\mathbf{x},{x}/\varepsilon)+ 
\varepsilon^2\theta^{(2)}(\mathbf{x},x/\varepsilon)+..., 
\label{theta2sc}
\\
\lambda^\varepsilon&\sim&\lambda_0.  
\label{lambda2sc}
\end{eqnarray}
Here 
$\mathbf{u}^{(j)}(\mathbf{x},y)$, $\theta^{(j)}(\mathbf{x},y)$, $j=0,1,2,...$, are functions 
to be determined of two independent variables: 
of a continuous macroscopic variable $\mathbf{x}$, and 
of the microscopic variable  $y=x/\varepsilon$ on the unit cell graph $\Gamma$.  All these functions are required to be 
one-periodic in $y_1$ and $y_2$ (which we henceforth conventionally call ``$\Gamma$-periodic in $y$'')  
 and 
to satisfy the usual kinematic compatibility conditions at the joints, cf. 
(\ref{contjunc}). 
When $\varepsilon\to 0$, the main-order terms $\mathbf{u}^{(0)}(\mathbf{x},y)$ and $\theta^{(0)}(\mathbf{x},y)$, 
are expected to describe the limit two-scale problem associated with $\lambda_0$, whose properties are of the 
main interest to us. 

The two-scale ansatz (\ref{u2sc})--(\ref{lambda2sc}) is formally substituted into (\ref{eqspvarmacro}) where 
the test functions are also chosen to be of a two-scale form: 
$\tilde{\mathbf{u}}^\varepsilon(\mathbf{x})\sim\tilde{\mathbf{u}}^{(0)}(\mathbf{x},{x}/\varepsilon)+
\varepsilon\tilde{\mathbf{u}}^{(1)}(\mathbf{x},{x}/\varepsilon)+..$ and 
$\tilde{\theta}^\varepsilon(\mathbf{x})\sim\tilde{\theta}^{(0)}(\mathbf{x},{x}/\varepsilon)+
\varepsilon\tilde{\theta}^{(1)}(\mathbf{x},{x}/\varepsilon)+..$, 
with $\tilde{\mathbf{u}}^{(j)}$ and $\tilde{\theta}^{(j)}$, $j=0,1,..$, $\Gamma$-periodic in $y$, satisfying kinematic compatibility at the joints in $y$, cf. (\ref{contjunc}), and having a bounded support 
in $\mathbf{x}$. 

Upon equating the terms of order $\varepsilon^{-2}$ in (\ref{eqspvarmacro}), and taking account of the 
arbitrariness of the above main-order ($j=0$) test functions, results in the following weak form condition for 
$\mathbf{u}^{(0)}(\mathbf{x},y)$ and $\theta^{(0)}(\mathbf{x},y)$. For any fixed $\mathbf{x}$, as functions of
$y$  both $\mathbf{u}^{(0)}(\mathbf{x},y)$ and $\theta^{(0)}(\mathbf{x},y)$ are a priori characterized by 
kinematically admissible $\Gamma$-periodic functions $(u(y), v(y),\theta(y))$  
 on the unit cell $\Gamma$, and the following integral 
identity holds for the stiff component $\Gamma_1$ of $\Gamma$ for all the admissible 
$\Gamma$-periodic test functions
 $(\tilde u(y), \tilde v(y),\tilde\theta(y))$: 
\begin{equation}
\int_{\Gamma_{1}}
\left[
\,\gamma\frac{du}{dy}\frac{d\tilde u}{dy}\,+\,\eta\frac{d\theta}{dy}\frac{d\tilde\theta}{dy}
\,+\, \kappa\left(\frac{d v}{dy}-\theta\right)
\left(\frac{d \tilde v}{dy}-\tilde\theta\right)
\right]\,dy\,\,=\, 0. 
\label{u0yeqn}
\end{equation}
Problem (\ref{u0yeqn}) plays a key role for determining the structure of the limit two-scale 
fields $\mathbf{u}^{(0)}(\mathbf{x},y)$ and $\theta^{(0)}(\mathbf{x},y)$, as follows. 

Setting in (\ref{u0yeqn}) $\left(\tilde u, \tilde v,\tilde\theta\right)=\left(\overline{u(y)}, \overline{v(y)},
\overline{\theta(y)}\right)$ (the overbar denotes complex conjugate) and 
as long as $\Gamma_{\rm stiff}$ is a {\it connected} periodic graph (as it is in the example of Fig. 1), 
(\ref{u0yeqn}) implies that $(u(y), v(y),\theta(y))$ can only correspond to a rigid body motion of 
of $\Gamma_1$, i.e. its translations and rotation with a constant $\theta$. 
The $\Gamma$-periodicity requirement further implies that on $\Gamma_1$,  
$\theta(y)\equiv 0$, and the total displacement 
$\mathbf{u}=u\mbox{\boldmath{$\tau$}}+v\mathbf{n}$, cf. (\ref{contjunc}), 
is independent of $y$ i.e. is a function of $\mathbf{x}$ only. On the soft part $\Gamma_0$, the functions 
$u(y)$, $v(y)$ and $\theta(y)$ remain at this point arbitrary. This allows to conclude that $\mathbf{u}^{(0)}$ 
and $\theta^{(0)}$ are of the form: 
\begin{eqnarray}
\mathbf{u}^{(0)}(\mathbf{x},y)&=&\mathbf{u}_0(\mathbf{x})\,+\,\mathbf{U}(\mathbf{x},y) 
\label{u0form} 
\\
\theta^{(0)}(\mathbf{x},y)&=& \ \ \ \ \ \ \ \ \ \ \ \  \,\Theta(\mathbf{x},y),  
\label{th0form}
\end{eqnarray} 
where $\mathbf{U}=(U,V)$ and $\Theta$ are supported only on the soft part $\Gamma_0$ 
(i.e. are identically zero on $\Gamma_1$) and satisfy zero (``clamping'') boundary 
conditions on the endpoints of $\Gamma_0$. 
The actual equations for the functions entering the right hand sides of (\ref{u0form}) and (\ref{th0form}) are 
still to be found, which are determined from equating next-order terms in the asymptotic procedure as 
follows. 

Equating next the terms of order $\varepsilon^{-1}$ in (\ref{eqspvarmacro}) leads to a corrector problem 
on the stiff part $\Gamma_1$ of the unit cell, as follows. For all kinematically admissible $\Gamma$-periodic test functions 
$\tilde u(y)$, $\tilde v(y)$ and $\tilde\theta(y)$, 
\begin{equation}
\int_{\Gamma_1}
\left[
\gamma\left(\frac{\partial\mathbf{u_0}}{\partial x}\cdot\mbox{\boldmath{$\tau$}}+
\eta\frac{du^{(1)}}{dy}\right)\frac{d\tilde u}{dy}+\eta\frac{d\theta^{(1)}}{dy}\frac{d\tilde\theta}{dy}
+ \kappa\left(\frac{\partial\mathbf{u_0}}{\partial x}\cdot\mathbf{n}+\frac{d v^{(1)}}{dy}-\theta^{(1)}\right)
\left(\frac{d \tilde v}{dy}-\tilde\theta\right)
\right]dy=0. 
\label{u1yeqn}
\end{equation}
Here, cf (\ref{eqspvarmacro}), $\partial\mathbf{u_0}/\partial x$ denotes the derivative of 
$\mathbf{u_0}(\mathbf{x})$ in the direction $dx$ along the beam; $\mbox{\boldmath{$\tau$}}$ is the unit vector in the direction 
of $dx$ and $\mathbf{n}$ is the unit vector in the normal direction of the positive transverse positive displacement $v$ on the 
beam; the unknowns $u^{(1)}$, $v^{(1)}$ and $\theta^{(1)}$ correspond to $\mathbf u^{(1)}$ for any fixed 
macroscopic variable $\mathbf{x}$ as functions on the unit cell graph $\Gamma$. 

With $\mathbf{u_0}=\left(u^0_1, u^0_2\right)$, $\mathbf{x}=(x_1,x_2)$, etc, and summation implied with 
respect to repeated indices, 
\[
\frac{\partial\mathbf{u_0}}{\partial x}\cdot\mbox{\boldmath{$\tau$}}\,=\,
\frac{\partial u^0_j}{\partial x_l}\tau_j\tau_l, \ \ \ 
\frac{\partial\mathbf{u_0}}{\partial x}\cdot\mathbf{n}\,=\,
\frac{\partial u^0_j}{\partial x_l}\tau_jn_l, 
\]
and so, from (\ref{u1yeqn}), $\mathbf{u}^{(1)}(\mathbf{x},y)$ and $\theta^{(1)}(\mathbf{x},y)$ 
can be represented in the form 
\begin{equation}
\left(\mathbf{u}^{(1)},\theta^{(1)}\right)(\mathbf{x},y)\,=\,
\frac{\partial u^0_j}{\partial x_l}(\mathbf{x})\mathbf{N}^{jl}(y), 
\label{cellcorr}
\end{equation} 
where $\mathbf{N}^{jl}(y)=\left(N^{jl}_u(y), N^{jl}_v(y), N^{jl}_\theta(y)\right)$, $j,l=1,2$, 
are solutions of the following unit-cell problems, in the weak form, 
\begin{equation}
\int_{\Gamma_{1}}
\left[
\gamma\left(\tau_j\tau_l+
\frac{d N^{jl}_u}{dy}\right)\frac{d\tilde u}{dy}+\eta\frac{d N^{jl}_\theta}{dy}\frac{d\tilde\theta}{dy}
+ \kappa\left(n_j\tau_l+
\frac{d N^{jl}_v}{dy}-N^{jl}_\theta\right)
\left(\frac{d \tilde v}{dy}-\tilde\theta\right)
\right]dy=0. 
\label{Nyeqn}
\end{equation}
For periodically connected stiff components $\Gamma_1$ the problem (\ref{Nyeqn}), which is a system 
of ordinary differential equations on $\Gamma_1$ with periodicity conditions determines, cf (\ref{u0yeqn}), 
$N^{jl}_\theta$ uniquely, and 
$N^{jl}_u$ and $N^{jl}_v$ up to arbitrary additive constants corresponding to rigid body translations, whose choice is insignificant. 
In simple cases like for the square network of Fig. \ref{fig3} it can be solved analytically, see Appendix. 
For purposes of the subsequent analysis, the correctors $\mathbf{N}^{jl}(y)$ should be regarded as 
as extended from $\Gamma_1$ to $\Gamma_0$ in an arbitrary kinematically compatible way, cf 
\cite{KSJMS18}, e.g. by linear interpolation for the soft segment on Fig. \ref{fig3}. 

Finally, to obtain the desired homogenized equations for $\mathbf{u}^{(0)}$, $\theta^{(0)}$, see 
(\ref{u0form})--(\ref{th0form}), it would suffice to take as test functions in (\ref{eqspvarmacro}) 
all those finitely supported in $\mathbf{x}$ of the form of (\ref{u0form})--(\ref{th0form}), (\ref{cellcorr}), i.e. 
\begin{eqnarray}
\tilde{\mathbf{u}}^\varepsilon(\mathbf{x})&=&\tilde{\mathbf{u}}_0(\mathbf{x})\,+\,
\varepsilon\tilde{\mathbf{u}}^{(1)}(\mathbf{x}, x/\varepsilon)
\,+\,\tilde{\mathbf{U}}(\mathbf{x},
x/\varepsilon) 
\label{utilform} 
\\
\tilde\theta^{\varepsilon}(\mathbf{x})&=& \varepsilon \theta^{(1)}(\mathbf{x}, x/\varepsilon)
\,+\,\tilde\Theta(\mathbf{x},x/\varepsilon),  
\label{thtilform}
\end{eqnarray} 
with the correctors 
\[
\left(\tilde{\mathbf{u}}^{(1)},\tilde\theta^{(1)}\right)(\mathbf{x},y)\,=\,
\frac{\partial \tilde u^0_p}{\partial x_q}(\mathbf{x})\left(N^{pq}_u(y), N^{pq}_v(y), N^{pq}_\theta(y)\right), 
\]
and $\tilde{\mathbf{U}}=(\tilde U,\tilde V)$ and $\tilde\Theta$ supported only on the soft part $\Gamma_0$ 
and vanishing on its border with $\Gamma_1$. 

Then (\ref{eqspvarmacro}), to the main order $\varepsilon^0$, results in 
\[
\int_{{\mathbb R}^2}C^h_{jlpq}\frac{\partial u^0_j}{\partial x_l}\frac{\partial \tilde u^0_p}{\partial x_q}dx\,+\,
\int_{{\mathbb R}^2}
\int_{\Gamma_{0}}
\left[
\,\frac{dU}{dy}\frac{d\tilde U}{dy}\,+\,\frac{d\Theta}{dy}\frac{d\tilde\Theta}{dy}
\,+\, \left(\frac{d V}{dy}-\Theta\right)
\left(\frac{d \tilde V}{dy}-\tilde\Theta\right)
\right]\,dy\,d\mathbf{x}\,=
\]
\begin{equation}
\lambda_0 \int_{{\mathbb R}^2}\int_{\Gamma}\left[ 
\left(\mathbf{u}_0 \cdot\mbox{\boldmath{$\tau$}}+U \right)
\left(\tilde{\mathbf{u}}_0 \cdot\mbox{\boldmath{$\tau$}}+\tilde U \right)+
\left(\mathbf{u}_0 \cdot\mathbf{n}+V \right)
\left(\tilde{\mathbf{u}}_0 \cdot\mathbf{n}+\tilde V \right)+
\Theta \tilde\Theta 
\right] \,dy\,d\mathbf{x}. 
\label{hom2sc}
\end{equation} 
In (\ref{hom2sc}), 
\[
C^h_{jlpq}\,:=\,\int_{\Gamma_1}\left[ \gamma 
\left(\tau_j\tau_l+\frac{d N^{jl}_u}{dy}\right)
\left(\tau_p\tau_q+\frac{d N^{pq}_u}{dy}\right)
\,+\,\eta\frac{d N^{jl}_\theta}{dy}\frac{d N^{pq}_\theta}{dy}\right.\ \ \ + 
\]
\begin{equation}
\ \ \ \ \ \ \ 
\left. 
+ \,\kappa\left(n_j\tau_l+\frac{d N^{jl}_v}{dy}-N^{jl}_\theta\right)
\left(n_p\tau_q+\frac{d N^{pq}_v}{dy}-N^{pq}_\theta\right)
\right]\,dy
\label{Chom}
\end{equation}
is the homogenized elasticity tensor for the stiff lattice $\Gamma_{\rm stiff}^\varepsilon$. 
One can show that, for periodically connected $\Gamma_1$, it is a conventional (generally anisotropic) elasticity tensor for a two-dimensional continuum 
linear elastic medium, satisfying all the usual conditions of symmetry and positivity, i.e. 
$C^h_{jlpq}=C^h_{pqjl}=C^h_{ljpq}$ and $E(e):=C^h_{jlpq}e_{jl}e_{pq}>0$ for any non-zero symmetric two-tensor 
$e$ (the strain tensor). 
In particular, for macroscopic rotations ($e_{lj}=-e_{jl}$) as $E(e)$ is the minimum of a non-negative quadratic 
functional corresponding to the right hand side of (\ref{Chom}), it is easily seen to vanish as minimized 
by $N_u\equiv N_v\equiv 0$ and $N_\theta=e_{jl}n_j\tau_l$ which is a microscopic rotation constant on 
$\Gamma_1$. 
Conversely, $E(e)=0$ implies $e_{lj}=-e_{jl}$, as the related minimizer can only correspond to a 
rigid body motion of $\Gamma_1$, and as can be confirmed by a direct inspection of (\ref{Chom}). 
Formula (\ref{Chom}), somewhat analogous to representations for homogenized tensors in classical 
periodic homogenization of continuum elastic media, expresses the homogenized tensor $C^h$ in terms of 
the solutions $N^{jl}_{u,v,\theta}$ of the cell problems (\ref{Nyeqn}). 
For simple geometries as e.g. the square lattice $C^h$ it 
can be computed analytically, see Appendix. 

The derived two-scale limit spectral problem (\ref{hom2sc}) is a system coupling the ``macroscopic'' part 
$\mathbf{u}_0(\mathbf{x})$ to 
the ``micro-resonant'' part 
$(\mathbf{U},\Theta)(\mathbf{x},y)=(U,V,\Theta)$, and 
is analogous to those originally derived by 
\cite{Zhikov00, Zhikov05} as adapted by us here to the network of high-contrast Timoshenko beams. We analyze some of its 
properties in the next section. 

\section{Analysis of the limit two-scale problem (\ref{hom2sc})}

Following the general recipe of \cite{Zhikov00,Zhikov05}, we aim at uncoupling (\ref{hom2sc}) by first expressing 
 $(\mathbf{U},\Theta)(\mathbf{x},y)=(U,V,\Theta)$ 
in terms of $\mathbf{u}_0(\mathbf{x})=\left(u^0_1,u^0_2\right)$. 
To this end, we set in (\ref{hom2sc}) $\tilde{\mathbf{u}}_0=0$ and notice that 
$\mathbf{u}_0\cdot\mbox{\boldmath{$\tau$}}=u^0_j(\mathbf{x})\tau_j$ and 
$\mathbf{u}_0\cdot\mathbf{n}=u^0_j(\mathbf{x})n_j$. This implies that 
\begin{equation}
U(\mathbf{x},y)=  u^0_j(\mathbf{x})U^j(y), \ \ \ 
V(\mathbf{x},y)=  u^0_j(\mathbf{x})V^j(y), \ \ \ 
\Theta(\mathbf{x},y)=  u^0_j(\mathbf{x})\Theta^j(y),
\label{vyform}
\end{equation}
where $\left(U^j, V^j, \Theta^j\right)$ solve microscopic problem on the soft part $\Gamma_0$ only: 
\[
\int_{\Gamma_{0}}
\left[
\,\frac{dU^j}{dy}\frac{d\tilde U}{dy}\,+\,\frac{d\Theta^j}{dy}\frac{d\tilde\Theta}{dy}
\,+\, \left(\frac{d V^j}{dy}-\Theta^j\right)
\left(\frac{d \tilde V}{dy}-\tilde\Theta\right)
\right]\,dy \,= 
\]
\begin{equation} 
\lambda_0  \int_{\Gamma_0}\left[ 
\left(\tau_j+U^j \right)\tilde U \,+ \,
\left(n_j+V^j \right)
\tilde V\,+\,
\Theta^j \tilde\Theta 
\right] \,dy. 
\label{vyeqn}
\end{equation} 
Equation (\ref{vyeqn}) forms a problem for real-valued $\left(U_j, V_j, \Theta^j\right)$ which are 
required to satisfy zero boundary conditions at the points of contact of $\Gamma_0$ with $\Gamma_1$ 
(the end points of the dashed segment on Fig. 3 in the example). This has a unique solution provided 
(real) $\lambda_0$ is not an eigenvalue 
of the (self-adjoint) spectral problem corresponding to (\ref{vyeqn}), i.e. if physically the applied frequency 
does not coincide with a resonant frequency of the soft phase.  
In this respect, equation (\ref{vyeqn}) can be viewed as explicitly accounting for the role of these micro-resonances. 

Returning to the limit two-scale problem (\ref{hom2sc}) and setting now 
$\tilde U=\tilde V=\tilde\Theta=0$, and recalling that for the test function $\tilde{\mathbf{u}}_0$, 
$\tilde{\mathbf{u}}_0\cdot\mbox{\boldmath{$\tau$}}=\tilde u^0_j(\mathbf{x})\tau_j$ and 
$\tilde{\mathbf{u}}_0\cdot\mathbf{n}=\tilde u^0_j(\mathbf{x})n_j$, we arrive at 
\begin{equation}
\int_{{\mathbb R}^2}C^h_{jlpq}\frac{\partial u^0_j}{\partial x_l}\frac{\partial \tilde u^0_p}{\partial x_q}dx\,
\,=\,
\lambda_0 \int_{{\mathbb R}^2}u^0_j(\mathbf{x})\tilde u^0_p(\mathbf{x})\int_{\Gamma}\left[ 
\left(\tau_j+U^j \right)\tau_p\,+\,
\left(n_j+V^j \right)n_p 
\right] \,dy\,d\mathbf{x}. 
\label{u0hom}
\end{equation} 
Integrating by parts, the latter transforms into the following partial differential equation for 
$\mathbf{u}_0(\mathbf{x})=(u_1,u_2)$: 
\begin{equation}
C^h_{jlpq}\frac{\partial^2u^0_p}{\partial x_l\partial x_q}\,\,+\,\, 
\beta_{jp}(\lambda_0)\,u^0_p\,\,=\,\,0, 
\label{u0pde}
\end{equation} 
where 
\begin{equation}
\beta_{jp}(\lambda)\,\,:=\,\,\lambda\,\int_{\Gamma}\left[ 
\left(\tau_j+U^j \right)\tau_p\,+\,
\left(n_j+V^j \right)n_p   
\right] \,dy. 
\label{beta}
\end{equation} 
In (\ref{beta}) the integration is performed over the whole of the unit cell graph $\Gamma$, with 
$U^p$ and $V^p$ extended by zero outside $\Gamma_0$ i.e. on the stiff component $\Gamma_1$. 

Setting in (\ref{vyeqn}) $(\tilde U,\tilde V,\tilde\Theta)=(U^p,V^p,\Theta^p)$, 
one observes that 
the above $2\times 2$ real-valued matrix $\beta(\lambda)$ is symmetric. Crucially, its signature for a 
given $\lambda$ determines whether the related frequency is in a propagating band or in a forbidden 
band gap, as follows. 

We seek a plane wave solution to (\ref{u0pde}), 
\[
\mathbf{u}_0(\mathbf{x})\,=\,\mathbf{A}e^{i\mathbf{k}\cdot\mathbf{x}}, 
\]
with a non-zero vector amplitude $\mathbf{A}=\left(A_1,A_2\right)$ and a wave vector $\mathbf{k}=(k_1,k_2)$. 
Substituting into (\ref{u0pde}) results in the following dispersion relation: 
\begin{equation}
C^h_{jlpq}k_lk_qA_p\,\,=\,\,\beta_{jp}(\lambda_0)A_p. 
\label{disprel}
\end{equation} 

As $C^h$ is positive definite, it immediately follows from (\ref{disprel}) that as long as the matrix 
$\beta(\lambda_0)$ is negative definite (i.e. both of its eigenvalues are negative) there is no 
non-trivial solution to (\ref{disprel}) and the corresponding frequency is in a bandgap. If on the other hand 
$\beta(\lambda_0)$ is positive definite, the related frequency is in a propagating band and (\ref{disprel}) 
implies the existence of two propagating modes in any direction as in a homogeneous linear elastic medium. 
An intermediate case of a ``weak gap'' occurs when $\beta(\lambda_0)$ is sign-indefinite i.e. when one of 
its eigenvalues is positive and the other is negative, cf \cite{Avila,VPS09,ZhP13}. In this case (\ref{disprel}) implies 
the existence of only one appropriately polarized propagating mode in any direction, with the medium 
thereby displaying some kind of polarization filtering effect. 

We remark finally that the above formal asymptotic derivations can in fact be stated and proved as a rigorous 
mathematical theorem, following the methodology of \cite{Zhikov00,Zhikov05}. Namely, for small enough 
$\varepsilon$, i.e. for high enough contrast, the above limit characterization of the bands and gaps 
describes the Floquet-Bloch spectrum of the original infinite periodic operator 
on $\Gamma^\varepsilon$ as determined by (\ref{eqspvarmacro}).

As we demonstrate in the next section, for simplest geometries including the one in our example of Figs. 1 and 
3, $\beta(\lambda)$ can be evaluate analytically. 

\section{Example: square network with a single isolated soft segment}

\begin{figure}
	\centering
		\includegraphics[scale=0.35]{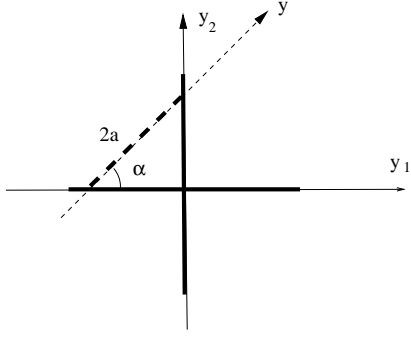}
		\caption{The unit cell graph for the example} 
		\label{fig4} 
\end{figure}

We specialize here the general results obtained in the previous section to the case of geometric configuration as 
in Figs. 1 and 3. Therefore the soft phase $\Gamma_0$ consists of an isolated inclined segment. 
Let the 
inclination angle be $\alpha$, $0^\circ<\alpha<90^\circ$, and the length of the segment be $2a$, Fig. \ref{fig4}. 

We start with evaluating the corresponding (symmetric) matrix $\beta(\lambda)$. It is convenient to 
equivalently evaluate quadratic form $\beta[\gamma]:=\beta_{jp}(\lambda)\gamma_j\gamma_p$ for an arbitrary vector 
$(\gamma_1,\gamma_2)$. 
 
In the present example, $\tau_1=n_2=\cos\alpha$, $\tau_2=-n_1=\sin\alpha$  are constant, 
and we choose the $y$-coordinate so that $\Gamma_0$ is the segment $[-a,a]$. 
Then (\ref{beta}) specializes to 
\begin{equation}
\beta[\gamma]=\,(2+2a)\lambda\left(\gamma_j\tau_j\right)^2+(2+2a)\lambda\left(\gamma_j n_j\right)^2\,+\,
\lambda\left(\gamma_p\tau_p\right)\int_{\Gamma_0}\gamma_jU^jdy\,+\, 
\lambda\left(\gamma_pn_p\right)\int_{\Gamma_0}\gamma_jV^jdy. 
\label{beta1}
\end{equation}
We further observe from (\ref{vyeqn}) that $\gamma_jU^j= \left(\gamma_j\tau_j\right)U$ 
and $\gamma_jV^j= \left(\gamma_j n_j\right)V$, where $U(y)$ and $V(y)$ are 
respectively the solutions of the following problems: 
\begin{equation}
\int_{-a}^a
\,\frac{dU}{dy}\frac{d\tilde U}{dy}\,dy \,= \, 
\lambda  \int_{-a}^a\,
\left(1\,+\,U \right)\tilde U \,dy,  
\label{Ueqn}
\end{equation} 
\begin{equation} 
\int_{-a}^a
 \,\left[
\frac{d\Theta}{dy}\frac{d\tilde\Theta}{dy}
\,+\, \left(\frac{d V}{dy}-\Theta\right)
\left(\frac{d \tilde V}{dy}-\tilde\Theta\right)\right]
\,dy \,= \, 
\lambda  \int_{-a}^a\left[ 
\left(1+V \right)
\tilde V\,+\,
\Theta \,\tilde\Theta 
\right] \,dy, 
\label{Veqn}
\end{equation} 
with zero boundary conditions 
$U(\pm a)=V(\pm a)=\Theta(\pm a)=\tilde U(\pm a)=\tilde V(\pm a)=\tilde\Theta(\pm a)=0$. 

As a result (\ref{beta1}) transforms into 
\begin{equation}
\beta[\gamma]\,=\,\lambda\left(\gamma_j\tau_j\right)^2\left[\,2\,+\,2a\,+\,
\int_{-a}^aU(y)\,dy\right]\,+\, 
\lambda\left(\gamma_jn_j\right)^2\left[\,2\,+\,2a\,+\,
\int_{-a}^aV(y)\,dy\right]. 
\label{beta2}
\end{equation}

Now, from (\ref{Ueqn}), $U(y)$ solves $U^{''}+\lambda(1+U)=0$ with boundary conditions 
$U(\pm a)=0$. As a result, as $\lambda\geq 0$, 
\begin{equation}
U(y)\,=\,\frac{\cos\lambda^{1/2}y}{\cos\lambda^{1/2}a}\,-\,1, 
\label{Usoln}
\end{equation}
provided $\lambda\neq\lambda_m:=\pi^2(m-1/2)^2/a^2$, $m=1,2,...$, i.e. is not an eigenvalue for longitudinal 
vibrations of the beam $\Gamma_0$. As a result, from (\ref{Usoln}), 
\begin{equation}
\int_{-a}^aU(y)\,dy\,=\, 2\lambda^{-1/2}\tan\lambda^{1/2}a\,-\,2a. 
\label{Umean}
\end{equation}

To evaluate $V(y)$, we observe from (\ref{Veqn}) that $V(y)$ and $\Theta(y)$ solve coupled system: 
\begin{eqnarray}
-\Theta^{''}\,-\,\left(V'-\Theta\right)&=&\lambda\Theta, \ \ \ 
V(\pm a)=\Theta(\pm a)=0, \nonumber\\ 
-\,\left(V'\,-\,\Theta\right)'&=&\lambda(1+V). 
\label{VTheqn}
\end{eqnarray}
Denoting $\hat V(y):=1+V(y)$ and eliminating  $\Theta$ in (\ref{VTheqn}), 
\begin{equation}
\hat V^{'''}\,+\,(1+\lambda)\hat V^{'}\,=\,(1-\lambda)\Theta,
\label{ThhatV}
\end{equation} 
\begin{equation}
\hat V^{{({\rm iv})}}\,+\,2\lambda\hat V^{''}\,+\,\lambda(\lambda-1)\hat V\,=\,0. 
\label{hatVeqn}
\end{equation} 
From the boundary conditions 
\begin{equation} 
\hat V(\pm a)\,=\,1,
\label{BC1}
\end{equation}
 and from (\ref{ThhatV}) together with $\Theta(\pm a)=0$, 
\begin{equation}
\hat V^{'''}(\pm a)\,+\,(1+\lambda)\hat V^{'}(\pm a)\,=\,0. 
\label{BC2}
\end{equation}

The characteristic equation for (\ref{hatVeqn}) is $\mu^4+2\lambda\mu^2+\lambda(\lambda-1)=0$, 
yielding $\mu^2=-\lambda\pm\lambda^{1/2}$. Assume first that $\lambda>1$, in which case 
$\mu=\pm i\left(\lambda\pm\lambda^{1/2}\right)^{1/2}$. The boundary conditions 
(\ref{BC1})--(\ref{BC2}) suggest that $V(y)$ is even, so is in the form 
\begin{equation}
\hat V(y)\,=\,A\cos\mu_1y\,+\,B\cos\mu_2y, \ \ \mu_1:=\left(\lambda+\lambda^{1/2}\right)^{1/2}, \ 
\mu_2:=\left(\lambda-\lambda^{1/2}\right)^{1/2}
\label{mu12}
\end{equation}
with constant $A$ and $B$. From (\ref{BC1}), 
\begin{equation}
A\cos\mu_1a\,+\,B\cos\mu_2a\,=\,1, 
\label{ABeq1}
\end{equation} 
and from (\ref{BC2}) 
\[
A\left(\mu_1^3-(\lambda+1)\mu_1\right)\sin\mu_1a\,+\,
B\left(\mu_2^3-(\lambda+1)\mu_2\right)\sin\mu_2a\,=\,0. 
\]
We further derive from the expressions (\ref{mu12}) for $\mu_1$ and $\mu_2$ that 
$\mu_1^3-(\lambda+1)\mu_1=\mu_1\mu_2^2\lambda^{-1/2}$ and 
$\mu_2^3-(\lambda+1)\mu_2 =-\mu_2\mu_1^2\lambda^{-1/2}$, 
and as a result, 
\begin{equation}
A\mu_2\sin\mu_1a\,-\,
B\mu_1\sin\mu_2a\,=\,0.
\label{ABeq2}
\end{equation} 
Provided $\mu_1\cos\mu_1a\sin\mu_2a+\mu_2\cos\mu_2a\sin\mu_1a\,\neq\,0$, i.e. the corresponding 
frequency is not an eigenfrequency of transverse vibrations of the clamped 
Timoshenko beam $\Gamma_0$, 
(\ref{ABeq1}) and (\ref{ABeq2}) determine $A$ and $B$ as follows: 
\[
A\,=\,\frac{\mu_1\sin\mu_2a}{\mu_1\cos\mu_1a\sin\mu_2a+\mu_2\cos\mu_2a\sin\mu_1a}, 
\]
\[
B\,=\,\frac{\mu_2\sin\mu_1a}{\mu_1\cos\mu_1a\sin\mu_2a+\mu_2\cos\mu_2a\sin\mu_1a}. 
\]
Finally, from (\ref{mu12}), 
\[
\int_{-a}^aV(y)dy\,=\,\int_{-a}^a\left(A\cos\mu_1y+B\cos\mu_2y\,-\,1\right)\,dy\,\,=\,
\]
\begin{equation}
\ \ \ \ \ \ \ \ \ \ \frac{4}{\mu_1\mbox{cotan}\,\mu_1a\,+\,\mu_2\mbox{cotan}\,\mu_2a}\,\,-\,\,2a.
\label{Vav}
\end{equation} 

As a result, substituting (\ref{Umean}) and (\ref{Vav}) into (\ref{beta2}) yields
\begin{equation}
\beta[\gamma]\,=\,\left(\gamma_j\tau_j\right)^2\beta_1(\lambda)\,\,+\, \, 
\left(\gamma_jn_j\right)^2\beta_2(\lambda), 
\label{betadiag}
\end{equation}
where 
\begin{equation}
\beta_1(\lambda)\,\,=\,\,2\lambda\,+\,2\lambda^{1/2}\tan\lambda^{1/2}a, 
\label{betau}
\end{equation} 
and 
\begin{equation}
\beta_2(\lambda)\,=\,2\lambda+
\frac{4\lambda}{\mu_1\mbox{cotan}\,\mu_1a+\mu_2\mbox{cotan}\,\mu_2a}, \ 
\mu_1:=\left(\lambda+\lambda^{1/2}\right)^{1/2}, \ 
\mu_2:=\left(\lambda-\lambda^{1/2}\right)^{1/2}. 
\label{betav}
\end{equation} 
Remark that for $0<\lambda<1$ the above calculation is still formally valid, with $\mu_2$ becoming 
imaginary and hence with relevant trigonometric functions replaced by their hyperbolic counterparts. 

We observe that $\beta_1(\lambda)$ and $\beta_2(\lambda)$ as given by (\ref{betau}) and (\ref{betav}) 
are precisely the two eigenvalues of the matrix $\beta(\lambda)$ as diagonalized by (\ref{betadiag}). 
As an important implication (see the discussion below (\ref{disprel})), the limit band gaps correspond 
to such frequencies that with associated $\lambda$ both $\beta_1(\lambda)$ and $\beta_2(\lambda)$ are 
negative. The existence of such gaps can be seen by direct inspection, i.e. by choosing the parameter $a$,  
and plotting $\beta_1(\lambda)$ and $\beta_2(\lambda)$ and detecting the domains of $\lambda$ where 
both are negative. 

The existence of such gaps can also be established qualitatively for large $\lambda$ and small $a$, 
by observing first 
that for $\lambda\gg 1$, asymptotically, 
\begin{equation}
\mu_{1,2}(\lambda)\,\sim\,\lambda^{1/2}\,\pm\,\frac{1}{2}\,-\,\frac{1}{8}\lambda^{-1/2}\,+\,O(\lambda^{-1}). 
\label{mu12asymp}
\end{equation}
Then, denoting $\Lambda:=\lambda^{1/2}$, 
from (\ref{betau}) we observe that $\beta_1(\lambda)$ is necessarily negative for 
$\Lambda=\pi/(2a)+(\pi/a)n+\delta_n$ for any positive integer $n$ and for small enough positive $\delta_n$. 
It is therefore sufficient to argue that, for appropriate $n$ and $a>0$, $\beta_2(\Lambda_n)<0$ where 
$\Lambda_n:=\pi/(2a)+(\pi/a)n$ are poles of $\beta_1(\Lambda)$. Using (\ref{mu12asymp}) and considering 
first large $n$ and hence large $\lambda=\lambda_n:=\Lambda_n^2$ one evaluates, asymptotically, 
\[
\mbox{cotan}(\mu_{1,2}a)\,\sim\,\mp\tan\frac{a}{2}\,+\,\frac{a}{8\Lambda_n}\tan^2\frac{a}{2}\,+\,
\frac{a}{8\Lambda_n}\,+\,
O\left(\Lambda_n^{-2}\right). 
\]
Then, substituting the above into (\ref{betav}), retaining main order terms in $\Lambda_n$, and 
finally assuming $a$ small results after straightforward calculations in 
\[
\beta_2(\Lambda_n)\,\sim\,2\Lambda_n^2\left(1\,-\,\frac{8}{a}\right). 
\]
So, for small enough $a$ and large enough $n$ both $\beta_1$ and $\beta_2$ are negative at 
$\Lambda_n+\delta_n$,  resulting in the band gaps. 

One can also easily detect the existence of ``weak'' gaps, when one of $\beta_j(\lambda)$, $j=1,2$,  
is positive and the other is negative. As a result, for the corresponding frequencies, in any particular direction only one appropriately 
polarized mode can propagate through such a medium, giving rise to some sort of polarization filtering. 
\vspace{.2in} 

{\bf Acknowledgments} 
\vspace{.15in}

The authors are thankful to the anonymous referee and to Prof J.R. Willis (University of Cambridge) for various comments and improving 
suggestions on the text. 
Thanks are also due to Prof A.S. Phani (University of British Columbia) and 
to Prof J. Kaplunov (Keele University) for useful discussions.

\appendix

\section{Homogenized elasticity tensor for the square network}

Here we calculate the homogenised elastic energy for the geometry considered in the previous section, 
Fig. \ref{fig4}.

Let us calculate $\mathbf{N}^{jl}:=\left({N}^{jl}_u, {N}^{jl}_v,{N}^{jl}_\theta\right)$. 
First let us notice that one can take $\mathbf{N}^{11}=0$ and  $\mathbf{N}^{22}=0$. Indeed,  for $\mathbf{N}^{11}$ (\ref{Nyeqn}) 
reduces to 
\begin{equation}
\int_{\Gamma_{1}}
\left[
\gamma\left(\tau_1\tau_1+
\frac{d N^{11}_u}{dy}\right)\frac{d\tilde u}{dy}+\eta\frac{d N^{11}_\theta}{dy}\frac{d\tilde\theta}{dy}
+ \kappa\left(
\frac{d N^{11}_v}{dy}-N^{11}_\theta\right)
\left(\frac{d \tilde v}{dy}-\tilde\theta\right)
\right]dy=0  
\label{IKNyeqn}
\end{equation}
since $\tau_1n_1=0$ on all the sides of the cross $\Gamma_1$. It remains to notice that 
$$
\int_{\Gamma_{1}}
\tau_1\tau_1
\frac{d\tilde u}{dy}=0
$$
since $\tau_1=0$ on $\Gamma_v$ (the vertical side of the cross) and $\tilde{u}$ is continuous and periodic on 
$\Gamma_h$ (the horizontal side). Consequently $\mathbf{N}^{11}=0$ is a solution. 
The same reasoning applies to $\mathbf{N}^{22}$.

Now we calculate $\mathbf{N}^{12}$. From (\ref{Nyeqn}) we have
\begin{equation}
\int_{\Gamma_{1}}
\left[
\gamma
\frac{d N^{12}_u}{dy}\frac{d\tilde u}{dy}+\eta\frac{d N^{12}_\theta}{dy}\frac{d\tilde\theta}{dy}
+ \kappa\left(n_1\tau_2+
\frac{d N^{12}_v}{dy}-N^{12}_\theta\right)
\left(\frac{d \tilde v}{dy}-\tilde\theta\right)
\right]dy=0. 
\label{IK2Nyeqn}
\end{equation}
By varying $\tilde u$, $N^{12}_u$ are constants on both $\Gamma_h$   and 
on $\Gamma_v$, which constants can be chosen both zero 
(as (\ref{IK2Nyeqn}) defines $\mathbf{N}^{12}$ up to a rigid translation). 
Next we notice that $n_1\tau_2=0$ on $\Gamma_{h}$  and $n_1\tau_2=-1$ on $\Gamma_{v}$  and (\ref{IK2Nyeqn}) reduces to 
$$
\int_{\Gamma_{h}}
\left[
\eta\frac{d N^{12}_\theta}{dy_1}\frac{d\tilde\theta}{dy_1}
+ \kappa\left(
\frac{d N^{12}_v}{dy_1}-N^{12}_\theta\right)
\left(\frac{d \tilde v}{dy_1}-\tilde\theta\right)
\right]dy_1\,+
$$
\begin{equation}
\int_{\Gamma_{v}}
\left[
\eta\frac{d N^{12}_\theta}{dy_2}\frac{d\tilde\theta}{dy_2}
+ \kappa\left(
\frac{d N^{12}_v}{dy_2}-N^{12}_\theta-1\right)
\left(\frac{d \tilde v}{dy_2}-\tilde\theta\right)
\right]dy_2=0.  
\label{IK3Nyeqn}
\end{equation}
We conclude (picking up zero  test function $\tilde \theta$  ) that 
\begin{equation}\frac{d N^{12}_v}{dy_1}-N^{12}_\theta=\alpha_h,\, \ \ \  \mbox{on}\, \ \ \Gamma_{h}\label{IK10}
\end{equation}
and 
\begin{equation}\frac{d N^{12}_v}{dy_2}-N^{12}_\theta=\alpha_v, \, \ \ \ \mbox{on}\, \ \ \Gamma_{v},\label{IK11}
\end{equation}
where $\alpha_h$ and $\alpha_v$ are some constants to be determined later.
Then the above formulas and (\ref{IK3Nyeqn}) (with $\tilde v=0$) imply
$$ \eta\frac{d^2N^{12}_\theta}{dy_1^2}\,=\,-\,\kappa\alpha_h,\,\ \ \  \mbox{on}\, \ \  \Gamma_{h}\setminus\{0\},
$$
and
$$ \eta\frac{d^2N^{12}_\theta}{dy_2^2}\,=\,\kappa(1-\,\alpha_v),\, \ \ \ \mbox{on}\, \ \ \Gamma_{v}\setminus\{0\}.
$$
The continuity and periodicity conditions for  $N^{12}_\theta$ imply
\begin{equation}  N^{12}_\theta=\frac{-\alpha_h}{2}\frac{\kappa}{\eta}\left(|y_1|-\frac{1}{2}\right)^2+\frac{\alpha_h}{8}\frac{\kappa}{\eta}  + B,\, \ \ \ \mbox{on}\, \ \ \Gamma_{h},
\label{IK8}\end{equation} 
and
\begin{equation}   N^{12}_\theta=\frac{\kappa}{\eta}\frac{(1-\alpha_v)}{2}\left(|y_2|-\frac{1}{2}\right)^2+
\frac{\kappa}{\eta}\frac{\alpha_v-1}{8}+B,\, \ \ \ \mbox{on}\, \ \ \Gamma_{v},
\label{IK9}\end{equation}  
where $B$ is some constant. 

Choosing $\tilde\theta=1$ and $\tilde{u}=0$ in (\ref{IK3Nyeqn}) we observe that 
\begin{equation}\alpha_h=1-\alpha_v\,.
\label{IK12}
\end{equation}
Next we note that (\ref{IK10}) implies
$$ \int_{\Gamma_{h}}\left(\alpha_h+N^{12}_\theta\right)
{dy_1}=0,
$$
since $N^{12}_v$ is periodic. Thus, using (\ref{IK8}) and integrating, 
\begin{equation}
0\,=\,\left(1+\frac{\kappa}{8\eta}\right)\alpha_h\,+\,B\,-\,
\frac{\kappa}{\eta}\int_{-1/2}^{1/2}\frac{\alpha_h}{2}\left(|y_1|-\frac{1}{2}\right)^2dy_1=
\left(1+\frac{\kappa}{12\eta}\right)\alpha_h+B.
\label{IK13}
\end{equation}
Similar considerations applied to (\ref{IK11}) imply, via (\ref{IK9}), 
$$0\,=\,\left(1+\frac{\kappa}{8\eta}\right)\alpha_v\,+\,B\,-\frac{\kappa}{8\eta}\,+\,
\frac{\kappa}{\eta}\int_{-1/2}^{1/2}\frac{1-\alpha_v}{2}\left(|y_2|-\frac{1}{2}\right)^2dy_2\,=
$$
\begin{equation}
\left(1+\frac{\kappa}{8\eta}\right)\alpha_v\,+\,B\,-\,\frac{\kappa}{8\eta}\,+\,
\frac{1}{24}(1-\alpha_v)\frac{\kappa}{\eta}\,=\,
\left(1+\frac{\kappa}{12\eta}\right)\alpha_v\,+\,B\,-\,\frac{\kappa}{12\eta}.
\label{IK14}
\end{equation}
Equations (\ref{IK12}), (\ref{IK13}) and (\ref{IK14}) imply
\begin{equation} \alpha_h=\frac{6\eta}{12\eta+\kappa},\, \ \ \ \alpha_v=\frac{6\eta+\kappa}{12\eta+\kappa},\, \ \ \ B=-\frac{1}{2}\,.
\label{IKconst}
\end{equation}
Thus $N^{12}_u$ and $N^{12}_\theta$ are fully determined. As for $N^{12}_v$,  one can find  an explicit expression using (\ref{IK10})  and (\ref{IK11}) but we will not use it in what follows.

Now we are ready to calculate $C^{h}_{jkpq}$ from its definition (\ref{Chom}). We have
$$C^h_{1111}=\int_{\Gamma_1}\left(\gamma\tau_1^2\tau_1^2+\kappa n_1\tau_1n_1 \tau_1\right) dy,
$$
since $\mathbf{N}^{11}=0$.  But $\tau_1n_1=0$; $\tau_1=0$ on $\Gamma_v$ and $\tau_1=1$ on $\Gamma_h$, therefore 

$$C^h_{1111}=\int_{\Gamma_{h}}\gamma dy\,=\,\gamma.
$$

Consider $C^h_{1112}$. We have
$$C^h_{1112}=\int_{\Gamma_1}\gamma\tau_1^2\tau_1\tau_2+
\kappa n_1\tau_1\left(n_1\tau_2+\frac{d N^{12}_v}{dy}-N^{12}_\theta\right)dy=0,
$$
since $N^{11}_u=N^{12}_u=N^{11}_\theta=N^{11}_v=0$.  But   $\tau_1\tau_2=\tau_1n_1=0$, consequently 
$$C^h_{1112}=0.
$$

Consider $C^h_{1122}$. We have
$$C^h_{1122}=\int_{\Gamma_1}\gamma\tau_1^2\tau_2^2+\kappa n_1\tau_1n_2\tau_2dy,
$$
since $\mathbf{N}^{11}=\mathbf{N}^{22}=0$. But  $\tau_1\tau_2=0$, therefore $C^h_{1122}=0$.

Consider $C^h_{1222}$. We have
$$C^h_{1222}=\int_{\Gamma_1}\gamma\tau_1\tau_2\tau_2^2+\kappa\left(n_1\tau_2+
\frac{d N^{12}_v}{dy}-N^{12}_\theta\right)
n_2\tau_2 dy,
$$
since $N^{12}_u=N^{22}_u=N^{22}_\theta=N^{22}_v=0$. 
But $\tau_1\tau_2=n_2\tau_2=0$ and consequently $C^h_{1222}=0$.

 For the ``shear'' term we have
$$C^h_{1212}=\int_{\Gamma_1}\left[ 
\eta\frac{d N^{12}_\theta}{dy}\frac{d N^{12}_\theta}{dy} +\kappa\left(n_1\tau_2+\frac{d N^{12}_v}{dy}-N^{12}_\theta\right)
\left(n_1\tau_2+\frac{d N^{12}_v}{dy}-N^{12}_\theta\right)
\right]\,dy=
$$
$$\int_{\Gamma_1} 
\kappa\left(n_1\tau_2+\frac{d N^{12}_v}{dy}-N^{12}_\theta\right)
n_1\tau_2
\,dy,
$$
where we have used (\ref{IK3Nyeqn}) with  $\tilde \theta= N^{12}_\theta$ and $\tilde v=N^{12}_v$. 
Since $n_1\tau_2=0$ on $\Gamma_h$ and $n_1\tau_2=-1$ on $\Gamma_v$ we have further simplification
$$C^h_{1212}=\int_{\Gamma_v}  
\kappa\left(-1+\frac{d N^{12}_v}{dy_2}-N^{12}_\theta\right)
\left(-1\right)
\,dy_2.
$$
Finally (\ref{IK11}) and (\ref{IKconst})  lead to
$$C^h_{1212}=\int_{\Gamma_v}  
\kappa\left(-1+\alpha_v\right)
\left(-1\right)\,dy=\frac{6\eta\kappa}{12\eta+\kappa}.
$$
For the last coefficient we have
$$C^h_{2222}=\int_{\Gamma_1}\left[ \gamma\tau_2^2\tau^2_2\,+\,\kappa\tau_2^2n_2^2\right]\,dy,
$$
since $\mathbf{N}^{22}=0$. Finally noticing that $\tau_2n_2=0$ on $\Gamma_1$, $\tau_2=0$ on $\Gamma_h$ and $\tau_2=1$ on $\Gamma_v$ we conclude
$$C^h_{2222}\,=\,\gamma.
$$
Substituting all the above evaluated values for $C^h_{jlpq}$ into 
the first term of (\ref{hom2sc}) we conclude that it equals to

\begin{equation}
\int_{{\mathbb R}^2}\left[\gamma\frac{\partial u^0_1}{\partial x_1}\frac{\partial \tilde u^0_1}{\partial x_1}+\gamma\frac{\partial u^0_2}{\partial x_2}\frac{\partial \tilde u^0_2}{\partial x_2}\,+\,
\frac{6\eta\kappa}{12\eta+\kappa}\left(\frac{\partial u^0_1}{\partial x_2}+\frac{\partial u^0_2}{\partial x_1}\right)
\left(\frac{\partial \tilde u^0_1}{\partial x_2}+\frac{\partial \tilde u^0_2}{\partial x_1}\right)\right]
dx.
\label{Chfin}
\end{equation}
Notice that, consistently with the original symmetries, the homogenized elasticity tensor $C^h$ as 
specified by (\ref{Chfin}) corresponds to an orthotropic two-dimensional elasticity tensor 
($C^h_{1111}=C^h_{2222}, \, C_{1112}=C_{2212}=0$), with zero Poisson-type ratio 
(i.e. $C_{1122}=0$).




\vspace{.2in} 


\end{document}